\documentclass[11pt]{amsart}

\usepackage{amsmath,amssymb}
\usepackage{microtype}
\usepackage[hidelinks]{hyperref}
\hypersetup{
 pdftitle={Strict inequalities for Stanley-Wilf limits},
 pdfauthor={Mikl\'os B\'ona}
}

\newcommand{\Av}{\operatorname{Av}}
\newcommand{\C}{\mathcal C}
\newcommand{\D}{\mathcal D}
\newtheorem{theorem}{Theorem}[section]
\newtheorem{lemma}[theorem]{Lemma}
\newtheorem{corollary}[theorem]{Corollary}
\theoremstyle{remark}
\newtheorem{remark}[theorem]{Remark}

\title[Growth rates of layered patterns]{Strict inequalities for Stanley-Wilf limits}
\author[M. B\'ona]{Mikl\'os B\'ona}
\address{Department of Mathematics, University of Florida,
Gainesville, FL 32611}
\email{bona@ufl.edu}
\date{}

\begin{document}

\begin{abstract}
Let $L(q)$ denote the Stanley--Wilf limit of a pattern $q$.
We prove that $L(1\oplus q)>(1+\sqrt{L(q)})^2$ for every skew
indecomposable pattern $q$ whose first entry is not $1$.
This strengthens the author's 2005 lower bound, which gives the
corresponding nonstrict inequality for every skew indecomposable
pattern. For $\delta_d=d(d-1)\cdots1$ and $d\ge2$, we obtain the
explicit bound $L(1\oplus\delta_d\oplus1)\ge\gamma_d>(d+1)^2$,
where $\gamma_d$ is the unique real number greater than $1$ satisfying
$\gamma_d^d(\gamma_d-1)=(d+1)^{2(d+1)}$.
In particular, \mbox{$L(14325)\ge16.2559832322\ldots>16$}.
We also prove that the Stanley--Wilf limit of every layered pattern
of length $n$ whose first and last layers have size one and which
has a layer of size greater than one exceeds $(n-1)^2$.
A refinement gives
$L(132546)>L(124356)>25$.
\end{abstract}

\maketitle

\section{Introduction}

A permutation $p$ contains a pattern $q$ if $p$ has a subsequence whose
entries occur in the same relative order as those of $q$. Otherwise,
we say that $p$ avoids $q$. Let $\Av_n(q)$ be the set of permutations
of length $n$ that avoid $q$, and let
\[
 a_n(q)=|\Av_n(q)|,
 \qquad
 L(q)=\lim_{n\to\infty}a_n(q)^{1/n}.
\]
The limit $L(q)$ exists and is finite for every fixed pattern
\cite{Arratia,MT}. It is called the Stanley--Wilf limit of $q$.
If $\iota_k=12\cdots k$, Regev's asymptotic formula
gives $L(\iota_k)=(k-1)^2$ \cite{Regev}.

Recall that the direct sum $\alpha\oplus\sigma$ places a copy of
$\alpha$ before and below a copy of $\sigma$. A permutation is
\emph{skew indecomposable} if it cannot be cut into two nonempty
parts so that every entry in the first part is larger than every
entry in the second part.

In 2005, the present author proved that
\[
 L(1\oplus q)\ge\left(1+\sqrt{L(q)}\right)^2
\]
for every skew indecomposable pattern $q$
\cite[Theorem~5.5]{Bona}. The most general result of this paper,
Theorem~\ref{thm:general-inflation}, shows that this inequality is
strict whenever the first entry of $q$ is not $1$:
\[
 L(1\oplus q)>\left(1+\sqrt{L(q)}\right)^2.
\]
When $q$ starts with $1$, equality holds by
\cite[Theorem~5.4]{Bona}. Thus Theorem~\ref{thm:general-inflation}
determines precisely when the 2005 lower bound is strict among
skew indecomposable patterns.

A permutation is called layered if it is a direct sum of decreasing
permutations. The same paper proved the lower bound
$L(q)\ge(k-1)^2$ for every layered pattern $q$ of length $k$
\cite{Bona}.

To introduce the construction, we first prove that this inequality
is strict for the layered patterns with layer sizes $1,d,1$, where
$d\ge2$. Let
\[
 \delta_d=d(d-1)\cdots1,
 \qquad
 \beta_d=1\oplus\delta_d\oplus1.
\]
Thus $\beta_2=1324$ and $\beta_3=14325$.
Our first result is the following.

\begin{theorem}\label{thm:main}
For every integer $d\ge2$, let $\gamma_d$ be the unique real number
greater than $1$ satisfying
\begin{equation}\label{eq:gamma}
 \gamma_d^d(\gamma_d-1)=(d+1)^{2(d+1)}.
\end{equation}
Then
\begin{equation}\label{eq:main}
 L(\beta_d)\ge\gamma_d>(d+1)^2=L(\iota_{d+2}).
\end{equation}
\end{theorem}

\begin{samepage}
\begin{corollary}\label{cor:14325}
Let $\gamma>1$ be the unique positive real solution of
\[
 \gamma^4-\gamma^3-65536=0.
\]
Then
\[
 L(14325)\ge\gamma=16.2559832322868\ldots>16.
\]
\end{corollary}
\end{samepage}

We will use the following known equalities:
\begin{equation}\label{eq:Dgrowth}
 |\Av_n(\delta_d\oplus1)|=|\Av_n(\iota_{d+1})|,
 \qquad
 L(\delta_d\oplus1)=d^2.
\end{equation}
The first assertion is the initial-block reversal equivalence of
Backelin, West, and Xin \cite{BWX}; the second follows from Regev's
formula. The next three sections give the proof.

In Section~\ref{sec:minima}, we extend the construction by showing
that increasing inflation at left-to-right minima preserves
avoidance under a weaker hypothesis than inflation at arbitrary
entries. This yields Theorem~\ref{thm:general-inflation} and the
strict improvement of the 2005 bound described above.
Section~\ref{sec:layered} applies this theorem to prove
$L(q)>(n-1)^2$ for every layered pattern $q$
of length $n$ with singleton first and last layers and at least one
layer of size greater than one. Later singleton layers are allowed.
We also obtain a strict recurrence for the family
$1\oplus21\oplus\cdots\oplus21\oplus1$.

\section{A decreasing staircase}

Let $d\ge2$ be fixed, and let
\[
 \C=\Av(12),
 \qquad
 \D=\Av(\delta_d\oplus1).
\]
We first construct permutations using a decreasing staircase of cells.
For an integer $r\ge1$, we use $r$ columns and $r+1$ rows, with rows
numbered $0,1,\ldots,r$ from top to bottom. In column $i$, we put a
$\D$-cell $D_i$ in row $i-1$ and a $\C$-cell $C_i$ in row $i$.
All other cells are empty. Figure~\ref{fig:staircase} illustrates the
case $r=4$.

\begin{figure}[ht]
\centering
\renewcommand{\arraystretch}{1.15}
\setlength{\arraycolsep}{14pt}
\[
\begin{array}{cccc}
 D_1 & & &\\
 C_1 & D_2 & &\\
 & C_2 & D_3 &\\
 & & C_3 & D_4\\
 & & & C_4
\end{array}
\]
\caption{The staircase for $r=4$. Higher rows contain larger values.
Each $C_i$ is decreasing, and each $D_i$ avoids $\delta_d\oplus1$.}
\label{fig:staircase}
\end{figure}

A gridded permutation consists of a permutation together with its
horizontal and vertical cuts. Positions can interleave arbitrarily
between the two cells in a column, and values can interleave arbitrarily
between the two cells in a row. Thus $C_i$ is directly below $D_i$,
and, when $i<r$, directly left of $D_{i+1}$.

\begin{lemma}\label{lem:avoidance}
Every permutation admitting this staircase gridding avoids $\beta_d$.
\end{lemma}

\begin{proof}
Suppose that an occurrence of $\beta_d$ exists. Its first entry is its
minimum and its last entry is its maximum. All its points therefore
lie in the rectangle between these two endpoints.

Two cells sharing neither a row nor a column are arranged in decreasing
order: the cell farther left is higher. An increasing pair consequently
lies either in one cell or in two cells sharing a row or a column.
The endpoint rectangle is contained in that cell or pair of cells.
Hence the entire occurrence lies in one cell or an adjacent pair.

A single $C_i$ cannot contain the occurrence because it is decreasing.
A single $D_i$ cannot contain it because $\beta_d$ contains
$\delta_d\oplus1$.

Consider a horizontal pair $C_i,D_{i+1}$. If at least two occurrence
points lie in $C_i$, they include the first two entries of $\beta_d$,
which form an ascent. Otherwise, at least the last $d+1$ entries lie
in $D_{i+1}$ and form $\delta_d\oplus1$.

For a vertical pair $C_i,D_i$, if at least two occurrence points lie
in the lower cell $C_i$, they include the two smallest values of
$\beta_d$. Those two values occur in increasing positional order.
Otherwise, the $d+1$ largest occurrence values lie in $D_i$ and form
$\delta_d\oplus1$. Every case is impossible.
\end{proof}

We now count permutations with a specified gridding of this form. Let
\[
 q_n^{(d)}=|\Av_n(\delta_d\oplus1)|.
\]
Let each $C_i$ contain exactly $m$ entries, and let each $D_i$ contain
exactly $dm$ entries.
There is one choice of permutation in each $C_i$ and
$q_{dm}^{(d)}$ choices in each $D_i$.
Each of the $r$ shared columns allows
$\binom{(d+1)m}{m}$ position interleavings, and each of the $r-1$
shared rows allows the same number of value interleavings.
These choices are independent. Indeed, once the positions and values
belonging to a cell have been chosen, the permutation in that cell
determines how to pair them. The number of
gridded base permutations is therefore
\begin{equation}\label{eq:B}
 B_{r,m}^{(d)}
 =\bigl(q_{dm}^{(d)}\bigr)^r
   \binom{(d+1)m}{m}^{2r-1}.
\end{equation}
Their length and the number of entries in the $C$-cells are
\begin{equation}\label{eq:NK}
 N=(d+1)rm,
 \qquad K=rm=\frac{N}{d+1}.
\end{equation}

Let
\[
 Q_d=\frac{(d+1)^{d+1}}{d^d}.
\]
By \eqref{eq:Dgrowth} and Stirling's formula, for each fixed $r$,
\begin{align}
 b_{d,r}
 &:=\lim_{m\to\infty}\bigl(B_{r,m}^{(d)}\bigr)^{1/((d+1)rm)}
 \notag\\
 &=d^{2d/(d+1)}Q_d^{(2r-1)/((d+1)r)}
  =(d+1)^2Q_d^{-1/((d+1)r)}.
 \label{eq:base-growth}
\end{align}
In particular, $b_{d,r}\to(d+1)^2$ as $r\to\infty$.

\section{Increasing inflations}

An \emph{interval} of a permutation is a set of consecutive positions
whose values form a set of consecutive integers. To inflate an entry
by $\iota_\ell$ is to replace it with an increasing interval of length
$\ell$, retaining the relative order of all blocks in both coordinates.

\begin{lemma}\label{lem:inflation}
Suppose that a pattern $\tau$ has no increasing interval of length at
least two. Inflating any entries of a $\tau$-avoider by nonempty
increasing permutations preserves $\tau$-avoidance.
In particular, this holds for $\tau=\beta_d$ whenever $d\ge2$.
\end{lemma}

\begin{proof}
If an occurrence uses at least two points of one inflated block, all
the selected points from that block occupy consecutive positions and
consecutive ranks within the occurrence. They therefore form a
nonsingleton increasing interval of $\tau$, a contradiction.
Otherwise the occurrence uses at most one point per block, and
contracting the blocks gives an occurrence in the original permutation.

For $d\ge2$, the first ascent and the last ascent of $\beta_d$ both
have value difference $d$, and all its other adjacent pairs decrease.
Thus $\beta_d$ has no adjacent increasing pair of consecutive values,
and hence no nonsingleton increasing interval.
\end{proof}

Starting with a base permutation counted by \eqref{eq:B}, inflate each
of its $K$ entries in the $C$-cells by an arbitrary nonempty increasing
permutation. Leave the $D$-entries unchanged. Perform all position and
value interleavings \emph{before} these inflations, so every inflated
block is an interval of the entire resulting permutation. Assign each
block to the cell of its original entry. By
Lemmas~\ref{lem:avoidance} and~\ref{lem:inflation}, the output avoids
$\beta_d$.

\begin{lemma}\label{lem:recovery}
The resulting gridded permutation uniquely determines its base gridded
permutation and all inflation lengths.
\end{lemma}

\begin{proof}
The original entries in a $C_i$-cell decrease. Its output subsequence
is therefore a skew sum of increasing blocks: every entry of an
earlier block exceeds every entry of a later block. The blocks are
exactly the maximal increasing runs of the $C_i$-subsequence.
Contracting these runs in every $C_i$ and leaving all $D_i$ entries
unchanged recovers the base permutation and its gridding. The run
lengths recover the inflation lengths. This also recovers $m$, so
outputs belonging to different values of $m$ remain distinguishable
when their griddings are supplied.
\end{proof}

If the inflations add a total of $j$ entries, we distribute these
entries among the $K$ original entries of the $C$-cells. Each of the
$K$ entries may receive any nonnegative number of additional entries,
so the number of choices is $\binom{K+j-1}{j}$.
For a permutation of length $M$, the number of
possible staircase griddings is at most
\begin{equation}\label{eq:griddings}
 (M+1)^{2r-1},
\end{equation}
because there are $r-1$ vertical cuts and $r$ horizontal cuts, each
chosen from at most $M+1$ gaps. Consequently, for every $j\ge0$,
\begin{equation}\label{eq:finite-bound}
 a_{N+j}(\beta_d)
 \ge
 \frac{B_{r,m}^{(d)}\binom{K+j-1}{j}}
 {(N+j+1)^{2r-1}}.
\end{equation}

For fixed $d,r$, let $g_{r,M}^{(d)}$ count all gridded outputs of
length $M$, allowing every $m\ge1$. Lemma~\ref{lem:recovery} shows
that their generating function is exactly
\begin{equation}\label{eq:G}
 G_r^{(d)}(z)
 :=\sum_{M\ge0}g_{r,M}^{(d)}z^M
 =\sum_{m\ge1}B_{r,m}^{(d)}
    z^{(d+1)rm}(1-z)^{-rm}.
\end{equation}
Moreover, \eqref{eq:griddings} gives the coefficientwise bound
\begin{equation}\label{eq:coeff-bound}
 g_{r,M}^{(d)}\le(M+1)^{2r-1}a_M(\beta_d).
\end{equation}

\begin{proof}[Proof of Theorem~\ref{thm:main}]
Fix $d\ge2$ and write $L=L(\beta_d)$.
For completeness, $\beta_d$ is skew indecomposable because its first
entry is its minimum. Hence $\Av(\beta_d)$ is closed under skew sum,
so $a_{n+m}(\beta_d)\ge a_n(\beta_d)a_m(\beta_d)$. Fekete's lemma,
together with the Marcus--Tardos bound \cite{MT}, gives the finite
limit $L$.

First apply \eqref{eq:finite-bound} with $j=0$. For fixed $r$, take
$N$th roots and let $m\to\infty$. The polynomial denominator
disappears, and \eqref{eq:base-growth} gives $L\ge b_{d,r}$.
Letting $r\to\infty$ yields $L\ge(d+1)^2>1$.

For every real $z$ with $0<z<1/L$, the series
\[
 \sum_{M\ge0}(M+1)^{2r-1}a_M(\beta_d)z^M
\]
converges. By \eqref{eq:coeff-bound}, so does $G_r^{(d)}(z)$.
The $((d+1)rm)$th root of its $m$th summand tends, as $m\to\infty$,
to
\[
 b_{d,r}\,z(1-z)^{-1/(d+1)}.
\]
Convergence therefore implies
\[
 b_{d,r}\,z(1-z)^{-1/(d+1)}\le1.
\]
This holds for each fixed $r$. Now let $r\to\infty$ and use
\eqref{eq:base-growth} to obtain
\[
 (d+1)^2z(1-z)^{-1/(d+1)}\le1
 \qquad(0<z<1/L).
\]
Finally, let $z\uparrow1/L$. We have already proved that
$L\ge(d+1)^2>1$, so $1/L$ lies strictly between $0$ and $1$.
The elementary function
$z\mapsto(d+1)^2z(1-z)^{-1/(d+1)}$ is continuous on $(0,1)$,
and in particular at $z=1/L$. Since the preceding inequality holds
for every $0<z<1/L$, passing to the limit from below therefore gives
\[
 \frac{(d+1)^2}{L}
 \left(1-\frac1L\right)^{-1/(d+1)}\le1.
\]
This step uses only continuity of the elementary function; it does
not require convergence of the generating function at $z=1/L$.
Rearranging the last inequality gives
\begin{equation}\label{eq:L-bound}
 L^d(L-1)\ge(d+1)^{2(d+1)}.
\end{equation}

The function $x\mapsto x^d(x-1)$ is continuous and strictly increasing
on $(1,\infty)$, from $0$ to infinity. This proves the existence and
uniqueness of $\gamma_d$ in \eqref{eq:gamma}, and
\eqref{eq:L-bound} implies $L\ge\gamma_d$.
Writing $A=(d+1)^2$, we have
$A^d(A-1)<A^{d+1}=(d+1)^{2(d+1)}$.
Thus $\gamma_d>A$, as required.
\end{proof}

\begin{remark}
The order of limits is essential to the estimate: $r$ is fixed while
$m$ tends to infinity, and only afterward does $r$ tend to infinity.
Thus the gridding factor \eqref{eq:griddings} is polynomial at every
stage where exponential growth is extracted. Also, the increasing
blocks are not freely reinterleaved with neighboring cells after
inflation; they remain global intervals throughout the construction.
\end{remark}

\section{The patterns 1324 and 14325}

For $d=2$, the $D$-cells avoid $213$ and have growth $4$.
Using cell lengths in ratio $1:2$ gives
\[
 L(1324)\ge\gamma_2,
 \qquad
 \gamma_2^3-\gamma_2^2-729=0,
 \qquad
 \gamma_2=9.3459834077485\ldots>9.
\]
For $d=3$, the $D$-cells avoid $3214$ and have growth $9$.
The ratio $1:3$ gives Corollary~\ref{cor:14325}.
Thus the same construction produces strict improvements over both
monotone benchmarks. It extends a mechanism proving $L(1324)>9$;
it does not require a numerical estimate for the full value of
$L(1324)$.

We can also obtain an explicit lower bound for $L(14325)$ directly
from \eqref{eq:finite-bound}. Let $m$ run through multiples of
$15$ and choose $j=K/15$. Then $N=4K=60j$ and $N+j=61j$.
For each fixed $r$, take $(61j)$th roots in \eqref{eq:finite-bound}
and let $m\to\infty$. Since
\[
 \lim_{j\to\infty}\binom{16j-1}{j}^{1/j}
 =\frac{16^{16}}{15^{15}},
\]
we obtain
\[
 L(14325)
 \ge b_{3,r}^{60/61}
     \left(\frac{16^{16}}{15^{15}}\right)^{1/61}.
\]
Here the exponent $60/61$ is the ratio of the base length $N$ to
the inflated length $N+j$. Indeed, since $4rm=N=60j$,
\[
 \bigl(B_{r,m}^{(3)}\bigr)^{1/(61j)}
 =\left(\bigl(B_{r,m}^{(3)}\bigr)^{1/(4rm)}\right)^{60/61}
 \longrightarrow b_{3,r}^{60/61}.
\]
The gridding denominator contributes the factor
$(61j+1)^{(2r-1)/(61j)}$, which tends to $1$ because $r$ is fixed.

We now explain the subsequent limit in $r$. By
\eqref{eq:base-growth},
\[
 b_{3,r}=16\left(\frac{256}{27}\right)^{-1/(4r)}
 \longrightarrow16
 \qquad\text{as }r\to\infty.
\]
The preceding lower bound holds for every positive integer $r$,
while $L(14325)$ is independent of $r$. Thus $L(14325)$ is at least
the limit of the right-hand side. Letting $r\to\infty$ and
simplifying the powers of $16$ gives
\begin{equation}\label{eq:explicit}
 \begin{aligned}
 L(14325)
 &\ge16^{60/61}
       \left(\frac{16^{16}}{15^{15}}\right)^{1/61}\\
 &=\frac{16^{76/61}}{15^{15/61}}
  =16\left(\frac{16}{15}\right)^{15/61}\\
 &=16.2559476393\ldots>16.
 \end{aligned}
\end{equation}
The final inequality is strict because $16/15>1$ and $15/61>0$.
The generating-function argument retains all inflation lengths and
gives the slightly larger constant $\gamma_3$.

\section{Inflating left-to-right minima}\label{sec:minima}

Lemma~\ref{lem:inflation} allows inflation at arbitrary entries, but
requires that the forbidden pattern have no nonsingleton increasing
interval. In our construction, only entries in the $C$-cells are
inflated. These entries are left-to-right minima, which allows us to
weaken that requirement.

Recall that an entry of a permutation is a \emph{left-to-right minimum}
if it is smaller than every entry preceding it.

\begin{lemma}\label{lem:minima-inflation}
Let $\tau$ be a pattern with $\tau_1=1$ and $\tau_2\ne2$.
Inflating any chosen left-to-right minima of a $\tau$-avoiding
permutation by nonempty increasing permutations preserves
$\tau$-avoidance.
\end{lemma}

\begin{proof}
Suppose that the inflated permutation contains $\tau$, and that an
occurrence selects $t\ge2$ entries from one inflated block $I$.
Since the entry replaced by $I$ was a left-to-right minimum, every
entry preceding $I$ is larger than every entry of $I$. The occurrence
starts with its minimum, so it must begin inside $I$.

The block $I$ is an interval of positions. Therefore its $t$ selected
entries are the first $t$ entries of the occurrence. It is also an
interval of values. Every selected entry outside $I$ must lie above
the whole block, since an entry below it would be smaller than the
first entry of the occurrence. Consequently, the selected entries
of $I$ are the $t$ smallest entries of the occurrence. They increase,
so $\tau$ begins with $12\cdots t$, contrary to $\tau_2\ne2$.

Thus every occurrence uses at most one entry from each inflated
block. Contracting the blocks would give an occurrence in the
original permutation, a contradiction.
\end{proof}

Increasing intervals later in $\tau$ cause no difficulty in this
lemma. Only an initial increasing interval could be supplied by an
inflated left-to-right minimum.

A permutation is \emph{skew indecomposable} if it cannot be cut into
two nonempty parts so that every entry in the first part is larger
than every entry in the second part. We can now extend the staircase
argument to a general class of patterns.

\begin{theorem}\label{thm:general-inflation}
Let $\sigma$ be skew indecomposable, with first entry different from
$1$, and let $s\ge1$ be a real number satisfying $s^2\le L(\sigma)$.
Define $\gamma_s>1$ by
\begin{equation}\label{eq:gamma-real}
 \gamma_s^s(\gamma_s-1)=(s+1)^{2(s+1)}.
\end{equation}
This extends the notation in \eqref{eq:gamma} to real $s$.
Then
\begin{equation}\label{eq:general-bound}
 L(1\oplus\sigma)\ge\gamma_s>(s+1)^2.
\end{equation}
In particular,
\begin{equation}\label{eq:strict-extension}
 L(1\oplus\sigma)>
 \left(1+\sqrt{L(\sigma)}\right)^2.
\end{equation}
\end{theorem}

\begin{proof}
Put $\pi=1\oplus\sigma$ and use the staircase of
Figure~\ref{fig:staircase}, with decreasing $C$-cells and
$D$-cells belonging to $\Av(\sigma)$.
Every entry in a $C_i$-cell is a left-to-right minimum of the whole
permutation. Indeed, earlier columns lie above $C_i$, the other cell
in its column is above $C_i$, and earlier entries within $C_i$ are
larger because that cell is decreasing.

An occurrence of $\pi$ could use a $C$-entry only as its initial
minimum. Its remaining entries would form $\sigma$ among the
$D$-entries. The $D$-cells are arranged in skew order, so skew
indecomposability forces this entire copy of $\sigma$ into one
$D$-cell, where it is forbidden. Thus every base permutation avoids
$\pi$.

Now inflate every $C$-entry by a nonempty increasing permutation,
after all interleavings have been made. Since $\pi_1=1$ and
$\pi_2\ne2$, Lemma~\ref{lem:minima-inflation} preserves avoidance.
Lemma~\ref{lem:recovery} applies without change: its proof uses only
the fact that the original $C$-cells decrease. The gridding bound
\eqref{eq:griddings} is unchanged as well.

Fix $r$, put $u_m=\lfloor sm\rfloor$, and use $m$ entries in each
$C_i$ and $u_m$ entries in each $D_i$. The number of gridded base
permutations, their size, and their number of inflatable entries are
\begin{equation}\label{eq:general-count}
 A_{r,m}=a_{u_m}(\sigma)^r
   \binom{m+u_m}{m}^{2r-1},
 \qquad N_m=r(m+u_m),\qquad K_m=rm.
\end{equation}
Let
\[
 Q_s=\frac{(s+1)^{s+1}}{s^s},
 \qquad b_{s,r}=(s+1)^2Q_s^{-1/((s+1)r)}.
\]
Since $L(\sigma)\ge s^2$, Stirling's formula gives
\begin{equation}\label{eq:general-limits}
 \lim_{m\to\infty}A_{r,m}^{1/N_m}\ge b_{s,r},
 \qquad
 \frac{K_m}{N_m}\longrightarrow\frac1{s+1}.
\end{equation}
Here $u_m/m\to s$, so the integer parts do not affect the limits.
Also, the recovery argument still determines $m$ from the number
of increasing runs in a $C$-cell.

The generating function for all inflated gridded outputs is
\[
 F_r(z)=\sum_{m\ge1}A_{r,m}z^{N_m}(1-z)^{-K_m}.
\]
Writing $L=L(\pi)$, we have
\[
 [z^M]F_r(z)\le(M+1)^{2r-1}a_M(\pi).
\]
The uninflated count, followed by $r\to\infty$, first gives
$L\ge(s+1)^2>1$. For every $0<z<1/L$, the last coefficient bound
implies convergence of $F_r(z)$. Equation~\eqref{eq:general-limits}
therefore gives
\[
 b_{s,r}z(1-z)^{-1/(s+1)}\le1.
\]
For each fixed $z$ with $0<z<1/L$, this inequality holds for every
positive integer $r$. Since
$b_{s,r}=(s+1)^2Q_s^{-1/((s+1)r)}\to(s+1)^2$ as
$r\to\infty$, passing to the limit in $r$ gives
\[
 (s+1)^2z(1-z)^{-1/(s+1)}\le1
 \qquad(0<z<1/L).
\]
Now let $z\uparrow1/L$. We have already proved that
$L\ge(s+1)^2>1$, so $1/L$ lies strictly between $0$ and $1$.
The elementary function
$z\mapsto(s+1)^2z(1-z)^{-1/(s+1)}$ is continuous on $(0,1)$,
and in particular at $z=1/L$. The preceding inequality holds for
every $0<z<1/L$, so passing to the limit from below gives
\[
 \frac{(s+1)^2}{L}
 \left(1-\frac1L\right)^{-1/(s+1)}\le1.
\]
As in the proof of Theorem~\ref{thm:main}, this step uses only
continuity of the elementary function; it does not require
convergence of $F_r(z)$ at $z=1/L$. Raising both sides to the
positive power $s+1$ and rearranging yields
\[
 L^s(L-1)\ge(s+1)^{2(s+1)}.
\]
The function $x\mapsto x^s(x-1)$ is strictly increasing from $0$ to
infinity on $(1,\infty)$, also for real $s\ge1$. The same comparison
at $x=(s+1)^2$ as before proves \eqref{eq:general-bound}.
Finally choose $s=\sqrt{L(\sigma)}$ to obtain
\eqref{eq:strict-extension}.
\end{proof}

\section{Further layered patterns}\label{sec:layered}

We apply the preceding theorem to layered patterns whose internal
layers may now have arbitrary positive lengths.

\begin{theorem}\label{thm:layered-second}
Let $\pi$ be a layered pattern of length $n$ whose first and last
layers have size one and whose second layer has size greater than
one. Then
\begin{equation}\label{eq:layered-quantitative}
 L(\pi)\ge\gamma_{n-2}>(n-1)^2.
\end{equation}
\end{theorem}

\begin{proof}
Write $\pi=1\oplus\sigma$. If the second layer of $\pi$ has size
$a\ge2$, then $\pi_2=a+1$ and the first entry of $\sigma$ is $a$.
The last entry of $\sigma$ is its maximum, so $\sigma$ is skew
indecomposable. It is also layered and has length $n-1$.
The general lower bound for layered patterns \cite{Bona} gives
$L(\sigma)\ge(n-2)^2$.
Apply Theorem~\ref{thm:general-inflation} with $s=n-2$.
\end{proof}

For example, $13245768$ has layer sizes $1,2,1,1,2,1$.
It has the increasing interval $45$, but
Lemma~\ref{lem:minima-inflation} still applies. Thus
$L(13245768)\ge\gamma_6>49$.

The restriction on the second layer can be removed from the strict
comparison by using the known singleton-extension formula
\begin{equation}\label{eq:singleton-extension}
 \sqrt{L(1\oplus q)}=1+\sqrt{L(q)}
 \qquad\text{if the first entry of $q$ is $1$}
\end{equation}
from \cite{Bona}.

\begin{corollary}\label{cor:singleton-ends}
Let $\pi$ be a layered pattern of length $n$ whose first and last
layers have size one. If at least one layer has size greater than
one, then
\[
 L(\pi)>(n-1)^2.
\]
\end{corollary}

\begin{proof}
Suppose that $\pi$ has $h\ge1$ initial singleton layers before its
first layer of size greater than one. Remove the first $h-1$ entries
and standardize to obtain $\rho$. Its length is $n-h+1$, its first
and last layers are singletons, and its second layer has size
greater than one. Theorem~\ref{thm:layered-second} gives
$\sqrt{L(\rho)}>n-h$.
Restore the removed initial singleton layers one at a time. By
\eqref{eq:singleton-extension},
\[
 \sqrt{L(\pi)}=(h-1)+\sqrt{L(\rho)}>n-1.
\]
\end{proof}

For the family with layer sizes $1,2,\ldots,2,1$, we obtain a strict
recurrence.

\begin{corollary}\label{cor:two-layers}
For $k\ge1$, let
\[
 \tau_k=1\oplus
 \underbrace{21\oplus\cdots\oplus21}_{k-1\text{ copies}}
 \oplus1,
 \qquad L_k=L(\tau_k).
\]
Thus $\tau_1=12$, $\tau_2=1324$, and $\tau_3=132546$.
For every $k\ge2$,
\begin{equation}\label{eq:family-recurrence}
 L_k>\left(2+\sqrt{L_{k-1}}\right)^2.
\end{equation}
Consequently, $L_k>(2k-1)^2$ for every $k\ge2$.
\end{corollary}

\begin{proof}
Let $\sigma_k$ be obtained by deleting the initial minimum of
$\tau_k$ and standardizing. Reversing the initial block $21$ using
\cite{BWX} shows that $\sigma_k$ is Wilf-equivalent to
\[
 12\oplus
 \underbrace{21\oplus\cdots\oplus21}_{k-2\text{ copies}}
 \oplus1=1\oplus\tau_{k-1}.
\]
Therefore \eqref{eq:singleton-extension} gives
\[
 L(\sigma_k)=\left(1+\sqrt{L_{k-1}}\right)^2.
\]
The pattern $\sigma_k$ is skew indecomposable and starts with $2$.
Apply \eqref{eq:strict-extension} to obtain
\eqref{eq:family-recurrence}. Taking square roots and iterating from
$L_1=1$ proves the last assertion.
\end{proof}

In particular, the case $k=3$ gives
\begin{equation}\label{eq:132546}
 L(132546)>
 \left(2+\sqrt{L(1324)}\right)^2
 =L(124356)>25.
\end{equation}
For the equality, note that $124356=1\oplus1324\oplus1$ and apply
\eqref{eq:singleton-extension} at both ends; the terminal version
follows by reverse complementation. The last inequality already
follows from the bound for $L(1324)$ in Section~4.
Thus $a_n(132546)>a_n(124356)$ for all sufficiently large $n$,
with strictly larger exponential order.

\end{document}